\documentclass{article}
\usepackage{amssymb}
\newcommand{\eq}{\begin{equation}}
\newcommand{\en}{\end{equation}}
\newcommand{\eqs}{\begin{eqnarray*}}
\newcommand{\ens}{\end{eqnarray*}}
\newcommand{\eqa}{\begin{eqnarray}}
\newcommand{\ena}{\end{eqnarray}}

\newcommand{\ex}{\mathbb E}

\newtheorem{theorem}{\large Theorem}[section]
\newtheorem{proposition}[theorem]   {\large Proposition}

\newtheorem{corollary}[theorem]{\large  Corollary}
\newtheorem{lemma}[theorem]{\large  Lemma}

\def\numberlikeadb{\global\def\theequation{\thesection.\arabic{equation}}}
\numberlikeadb

\def\ignore#1{}
\def\E{{\mathbb E}}

\def\R{{\mathbb R}}
\def\real{\R}
\def\sli{\sum_{l\ge1}}
\def\cc{{\cal C}}
\def\Eq{\ =\ }
\def\Le{\ \le\ }
\def\Ref#1{(\ref{#1})}

\def\h{\eta}
\def\sji{\sum_{j\ge1}}
\def\dtv{d_{{\rm TV}}}
\def\Le{\ \le\ }
\def\law{{\cal L}}
\def\bone{{\bf 1}}
\def\nin{\noindent}

\def\pr{{\bf P}}
\def\half{{\textstyle{{1\over2}}}}

\def\b{\beta}
\def\e{\varepsilon}

\def\sjm{\sum_{j=1}^m}

\def\lti{\lim_{t\to\infty}}

\def\law{{\cal L}}

\def\mvn{{\rm MVN\,}}
\def\ep{\hfill$\Box$}
\def\nn{{\cal N}}
\def\var{{\rm Var\,}}
\def\cov{{\rm Cov\,}}

\def\a{\alpha}
\def\r{\rho}
\def\sii{\sum_{i\ge1}}
\def\f{\phi}
\def\sjo{\sum_{j\ge0}}
\def\Corr{{\rm Corr\,}}

\begin{document}

\title{Small counts in the infinite occupancy scheme}
\author{A.D. Barbour\thanks{Angewandte Mathematik,
Winterthurerstrasse~190,
CH--8057 Z\"urich, Switzerland: {\tt a.d.barbour@math.uzh.ch}
}~~ and~~ A.V. Gnedin\thanks{Department of Mathematics, Utrecht University,
PO Box 80010, 3508 TA Utrecht, The Netherlands: {\tt A.V.Gnedin@uu.nl}
}\\
{\it University of Z{\"u}rich and Utrecht University}}
\date{\today}
%{\empty}
\maketitle

\begin{abstract}
\noindent
The paper is concerned with the classical occupancy 
scheme with infinitely many boxes, in which~$n$ balls are thrown independently into  
boxes $1,2,\ldots$, with probability~$p_j$ of hitting the box~$j$, where
$p_1\geq p_2\geq\ldots>0$ and $\sum_{j=1}^\infty p_j=1$.
We establish joint normal approximation as $n\to\infty$
for the numbers of boxes containing $r_1,r_2,\ldots,r_m$ balls,
standardized in the natural way, assuming only that the variances of these 
counts all tend to infinity. The proof of this approximation is based on a
de-Poissonization lemma. We then review sufficient conditions for the
variances to tend to infinity.
Typically, the normal approximation does not mean convergence.
We show that the convergence of the full vector of $r$-counts 
only holds under a condition of regular variation, thus giving a complete 
characterization of possible limit correlation structures.

\end{abstract}

\section{Introduction}
\setcounter{equation}{0}
In the classical occupancy 
scheme with infinitely many boxes,  balls are thrown independently into  
boxes $1,2,\ldots$, with probability $p_j$ of hitting the box $j$, where
$p_1\geq p_2\geq\ldots>0$ and
 $\sum_{j=1}^\infty p_j=1$. The most studied quantity is the number of 
 boxes $K_n$ occupied by at least one out of the first~$n$ balls thrown.
It is known that for large $n$ the law of $K_n$ is asymptotically 
normal, provided that ${\rm Var}[K_n]\to\infty$; see
\cite{GHP, Hwang} for references and a survey of this and related results.
In this paper, we investigate the behaviour of the quantities $X_{n,r}$, 
the numbers of boxes hit by exactly $r$ out of the~$n$  balls, $r\ge1$.

Under a condition of regular variation, a multivariate CLT for the $X_{n,r}$'s 
was proved by Karlin \cite{Karlin}. Mikhailov \cite{Mikhailov} also studied 
the $X_{n,r}$'s, but in a situation where the $p_j$'s vary with $n$.
In this paper, we establish joint normal approximation as $n\to\infty$
for the variables $X_{n,r_1},\ldots, X_{n,r_m}$, centred and normalized,
assuming only that $\lim_{n\to\infty}\var X_{n,r_i} = \infty$ for each~$i$.
We also give examples to show that this condition is not enough to ensure
{\em convergence\/}, since the correlation matrices need not converge
as $n\to\infty$.  The asymptotic behaviour of the moments of the~$X_{n,r}$
is thus of key importance, and we discuss this under a number of simplifying
assumptions.

\par The behaviour of these moments, as also of those
of $K_n=\sum_{r=1}^\infty X_{n,r}$, 
depends on the way in which the frequencies $p_j$ decay to $0$.
In the case of power-like decay, $p_j\sim cj^{-1/\alpha}$ with $0<\alpha<1$,
it is known that, for each fixed~$k$, the moments $\E X_{n,r}^k$ 
have the same order of growth with~$n$ for every~$r$, 
and this is the same order of growth as that of $\E K_n^k$; moreover,
the limit distributions of $K_n$ and of $X_n:=(X_{n,1},X_{n,2},\ldots)$ 
are normal \cite{GHP,Karlin}. 
In contrast, for  a sequence of geometric frequencies $p_j=cq^j$ ($0<q<1$),
there is no way to scale the $X_{n,r}$'s to obtain a nontrivial limit 
distribution \cite{Louchard},
and the moments of $K_n$ have oscillatory asymptotics.
In a more general setting such that the $p_j$'s have exponential decay, 
the oscillatory behaviour of ${\rm Var}[K_n]$ is typical \cite{BGY}.
The spectrum of interesting possibilities is, however, much wider:   
for instance, frequencies  $p_j\sim ce^{-j^\beta}$, with $0<\beta<1$,
exhibit a decay intermediate between power and exponential.

Karlin's \cite{Karlin} multivariate CLT for $X_n$ applies when the index of regular 
variation is in the range $0<\alpha<1$.
We complement this by the analysis of the cases $\alpha=0$ and $\alpha=1$,
showing that for each $\alpha\in [0,1]$ there is exactly one possible normal limit.
 Finally, we prove that these one-parameter normal laws are the only possible limits of 
naturally scaled and centred $X_n$.
Specifically, we show that a regular variation condition holds  
if $\var X_{n,r}\to\infty$ for all $r$ and if all the correlations
$\{\Corr(X_{n,r},X_{n,s}),\,r,s\ge1\}$ converge.

%CLT holds with one of
%the normal limits determined above.

\section{Poissonization}
\setcounter{equation}{0}
 As in much previous work, we shall rely on a closely related occupancy scheme,
in which the balls are thrown into the boxes at the times of a unit Poisson
process. The advantage of this model is that, for every $t>0$, the 
processes $(N_j(t)\,,~ t\geq 0)$,  
counting the numbers of balls in boxes $j=1,2,\ldots$, are independent.
Let $Y_r(t)$ be the number of boxes occupied by exactly $r$ balls at time $t$.
In  view of the representation 
\begin{equation}\label{PoInd}
   Y_r(t)=\sum_{j=1}^\infty {\bf 1}[N_j(t)=r]
\end{equation}
with independent Bernoulli terms, it follows that
\eq\label{Y-normal}
   Y'_r(t)\ :=\ (Y_r(t)-{\mathbb E}[Y_r(t)])/\sqrt{{\rm Var}[Y_r(t)]}
    \ \to_d\ \nn(0,1) \quad\mbox{as}\quad  t\to\infty
\en
if and only if ${\rm Var}[Y_r(t)]\to\infty$. This suggests that normal
approximation can be approached most easily through the $Y_r(t)$, provided
that the de-Poissonization can be accomplished.  We now show that this is
indeed the case.

\par Let ${\cal L}(\cdot)$ denote the probability law of a random element,
$\dtv$ the distance in total variation.

\begin{lemma}\label{2.1} 
For any $m,k \in {\mathbb N}$ satisfying $m \le {1\over 2} np_k$,
we have
$$
   \dtv(\law(X_{n,1},\ldots,X_{n,m}),\law(Y_1(n),\ldots,Y_m(n))) 
	     \Le \pi_k + 2k e^{- np_k/10},
$$
where $\pi_j:=\sum_{i=j+1}^\infty p_i$.
\end{lemma}

\nin{\bf Proof.}\ We begin by noting that, in parallel to~\Ref{PoInd},
\eq\label{MNInd}
   X_{n,r}:=\sum_{j=1}^\infty {\bf 1}[M_{n,j}=r],
\en	 
where $M_{n,j}$ represents the number of balls out of the first~$n$ thrown
that fall into box~$j$.  Our proof uses lower truncation of the sums
\Ref{PoInd} and~\Ref{MNInd} that define $Y_r(n)$ and~$X_{n,r}$.

Since $M_{n,j} \sim {\rm Binomial}(n,p_j)$, it follows from the Chernoff
inequalities \cite{CL06} that, if
$m \le \half np_k$, then for $j\leq k$
\[
   {\mathbb P}[M_{n,j} \le m] \Le {\mathbb P}[M_{n,j} \le \half np_j]
	   \Le \exp\{- np_j/10\} \Le \exp\{- np_k/10\}, 
\]
since the~$p_j$ are decreasing, and $m \le \half np_k$;
and the same bound holds also for  $N_j(n)\sim {\rm Poisson}(np_j)$.  
Hence, defining
$$
   X_{n,k,r} := \sum_{j=k+1}^\infty \bone[M_{n,j}=r],\qquad 
	 Y_{k,r}(t) := \sum_{j=k+1}^\infty\bone[N_j(t)=r],
$$
it follows that
\begin{eqnarray}
    \dtv(\law(X_{n,1},\ldots,X_{n,m}),\law(X_{n,k,1},\ldots,X_{n,k,m}))
		     & \le & ke^{- np_k/10}; \label{1}\\ 
    \dtv(\law(Y_1(t),\ldots,Y_m(t)),\law(Y_{k,1}(t),\ldots,Y_{k,m}(t)))
		     & \le & ke^{- tp_k/10}.\label{2}
\end{eqnarray}
But now, from an
inequality of Le Cam \cite{LeCam} and Michel \cite{Michel}, we have
\begin{equation}\label{3}
   \dtv(\law(N_j(n),\,j\ge k+1), \law(M_{n,j},\,j\ge k+1))  \Le \pi_k, 
\end{equation}
and the $X_{n,k,r}$ are functions of $\{M_{n,j},\,j\ge k+1\}$, the $Y_{k,r}(n)$
of $\{N_j(n),\,j\ge k+1\}$. The lemma now follows from  (\ref{1}),(\ref{2}) and (\ref{3}).
\ep

\bigskip
\noindent 

\begin{proposition}\label{2.2} 
Let $k(n)$ be any sequence satisfying 
$$
  k(n) \to \infty~~~{\rm and}~~~ k(n)e^{- np_{k(n)}/10} \to 0.
$$  
Then, for any sequence $m(n)$ satisfying
$m(n) \le \half np_{k(n)}$ for each~$n$, it follows that
\begin{equation}\label{4}
   \dtv(\law(X_{n,1},\ldots,X_{n,m(n)}),\law(Y_1(n),\ldots,Y_{m(n)}(n))) \ \to\ 0.
\end{equation}  
\end{proposition}

\nin{\bf Proof.}\ Since $m(n) \le \half np_{k(n)}$ for each~$n$,
it follows that Lemma~\ref{2.1} can be applied for each~$n$.  
Since $k(n) \to \infty$, 
it follows that $\pi_{k(n)} \to 0$, so that the first element in its bound
converges to zero; the second converges to zero also, by assumption.
\ep

\bigskip

\nin{\bf Remark.}\ Such sequences~$k(n)$ always exist. For instance, one can
take 
$$
   k(n) = \max\{k\colon\ 20\log k/p_k \le n\}.
$$
For this choice, it is immediate that ${k(n)} \to \infty$, and that
$np_{k(n)} \ge 20\log k(n) \to \infty$, entailing also that
$k(n)e^{- np_{k(n)}/10} \le 1/k(n) \to 0$. Hence there are always sequences
$m(n) \to \infty$ for which (\ref{4}) is satisfied.
\bigskip

\par Hence, in particular, any approximation to the distribution of a finite subset 
of the components of~$Y(n)=(Y_1(n),Y_2(n),\ldots)$ (suitably scaled) remains valid for the corresponding
components of~$X_n$, at the cost of introducing an extra, asymptotically
negligible, error in total variation of at most
\eq\label{tv-error}
   \pi_{k(n)} + 2k(n)e^{- np_{k(n)}/10},
\en
where $k(n)$ is any sequence satisfying the conditions of Proposition~\ref{2.2}.

\section{Normal approximation}\label{normal}
\setcounter{equation}{0}
As noted above, the distribution of $Y_r(t)$ is asymptotically normal 
as $t\to\infty$ whenever $\var Y_r(t) \to \infty$.  Here, we consider the joint normal
approximation of any finite set of counts $Y_{r_1}(t),\ldots, Y_{r_m}(t)$ such
that $r_i\ge1$ and $\lti\var Y_{r_i}(t) = \infty$ for each $1\le i\le m$.  We measure the closeness
of two probability measures $P$ and~$Q$ on~$\real^m$ in terms of differences between 
the probabilities assigned to arbitrary convex sets:
\[
     d_c(P,Q) \ :=\ \sup_{A\in\cc}|P(A)-Q(A)|,
\]
where $\cc$ denotes the class of convex subsets of~$\real^m$.
Let
\[
   \Phi_r(t) \ :=\ {\mathbb E} Y_r(t), ~~~V_r(t)\ :=\ \var Y_r(t),
    \quad C_{rs}(t)\ :=\ \cov(Y_r(t),Y_{s}(t)) 
\] 
denote the moments of the~$Y_r(t)$, and let 
\[
    \Sigma_{rs}(t)\ :=\ C_{rs}(t)/\sqrt{V_r(t)V_s(t)}\ =\ \cov(Y'_r(t),Y'_s(t))
\]
denote the covariance matrix of the standardized random variables~$Y'_r(t)$  
as in~\Ref{Y-normal}. 

Now the random vector
$(Y'_{r_1}(t),\ldots,Y'_{r_m}(t))$ is a sum of independent
mean zero random vectors $(Y'_{l,r_1}(t),\ldots,Y'_{l,r_m}(t))$, $l\ge1$,
where $Y'_{l,r}(t) := ({\bf 1}[N_l(t)=r] - p_{l,r}(t))/\sqrt{V_r(t)}$,
and 
\begin{equation}\label{plr}
p_{l,r}(t) := {\mathbb P}[N_l(t)=r]=e^{-tp_l}{(tp_l)^r\over r!}.
\end{equation}
A theorem of Bentkus \cite[Thm.~1.1]{B04} then shows that
\[
   d_c(\law(Y'_{r_1}(t),\ldots,Y'_{r_m}(t)),\mvn_m(0,\Sigma_R(t))) \Le
      C m^{1/4}\b_t,
\]
for an absolute constant~$C$, where
\[
   \b_t \ :=\ \sli \b_{t,l} \quad\mbox{and}\quad
    \b_{t,l} \ :=\ \ex|\Sigma^{-1/2}_R(t) (Y'_{l,r_1}(t),\ldots,Y'_{l,r_m}(t))^T|^3,
\]
and $\Sigma_R(t)$ denotes the $m\times m$ matrix with elements $\{\Sigma_{rs}(t),\,
r,s \in R := \{r_1,\ldots,r_m\}\}$.  Applying this result, we obtain the following theorem. 

\begin{theorem}\label{approximation}
If\, $\lti V_{r_i}(t) = \infty$ for each $1\le i\le m$, where $1 \le r_1 < \ldots < r_m$,
then, as $t$ and~$n$ tend to $\infty$,
\eqs
   d_c(\law(Y'_{r_1}(t),\ldots,Y'_{r_m}(t)),\mvn_m(0,\Sigma_R(t)))
     &=& O\Bigl(1\Big/\min_{1\le i\le m} \sqrt{V_{r_i}(t)}\Bigr) \ \to\ 0;\\
   d_c(\law(X'_{n,r_1},\ldots,X'_{n,r_m}),\mvn_m(0,\Sigma_R(n)))
     &=& O\Bigl(\pi_{k(n)} + 2k(n)e^{- np_{k(n)}/10} +
        \Bigl\{1\Big/\min_{1\le i\le m} \sqrt{V_{r_i}(n)}\Bigr\} \Bigr) \\
     &\to& 0,
\ens
where~$k(n)$ is any sequence chosen as for {\rm \, Proposition~\ref{2.2}} and satisfying
$\max_{1\le j\le m}r_j \le \half np_{k(n)}$ for each~$n$.
If, in addition, $\Sigma_R(t) \to \Sigma_R$ as $t\to\infty$, for some fixed~$\Sigma_R$,
then 
$$
   (Y'_{r_1}(t),\ldots,Y'_{r_m}(t))\ \to_d\ \mvn_m(0,\Sigma_R) \quad \mbox{and}
   \quad (X'_{n,r_1},\ldots,X'_{n,r_m})\ \to_d\ \mvn_m(0,\Sigma_R).
$$
\end{theorem}

\nin{\bf Proof.}\ All that we need to do is to control the quantity~$\b_t$. This
in turn involves bounding the smallest eigenvalue of~$\Sigma_R(t)$ away from~$0$.
Now direct calculation shows that, for any column vector $a \in \real^m$,
\[
   a^T \Sigma_R(t)\, a \Eq \var\Bigl(\sjm a_j Y'_{r_j}(t)\Bigr)
    \Eq \sli\var\Bigl(\sjm a_j Y'_{l,r_j}(t)\Bigr).
\]
Using the definition of~$Y'_{l,r}(t)$, this gives
\eqs
   a^T \Sigma_R(t)\, a &=& \sli\Bigl\{p^{l,R}(t)\ex^{l,R,t}(U^2) 
     - \{p^{l,R}(t)\}^2\,\{\ex^{l,R,t}(U)\}^2\Bigr\},
\ens
where $p^{l,R}(t) := \sum_{r\in R}p_{l,r}(t)$ and, under the measure~$\pr^{l,R,t}$, 
$U$ takes the value $a_j/\sqrt{V_{r_j}(t)}$ with probability $p_{l,r_j}(t)/p^{l,R}(t)$,
$1\le j\le m$.  This in turn implies that
\eqs
     a^T \Sigma_R(t)\, a &\ge& \sli p^{l,R}(t)(1-p^{l,R}(t))\ex^{l,R,t}(U^2),
\ens
and since 
\[
    \ex^{l,R,t}(U^2) = \sjm \frac{p_{l,r_j}(t)\, a_j^2}{p^{l,R}(t) V_{r_j}(t)},
\]
it follows that
\[
   a^T \Sigma_R(t)\, a \ \ge\ \sli (1-p^{l,R}(t))\sjm \frac{p_{l,r_j}(t)\, a_j^2}{V_{r_j}(t)}
    \ \ge\ \min_{l\ge1}(1-p^{l,R}(t))\, a^Ta,
\]
since $V_r(t) \le \sli p_{l,r}(t)$.  However, for each~$l$,
$p^{l,R}(t) \le 1 - p_{l,0}(t) - \sum_{j > r_m} p_{l,j}(t)$, and 
$(p_{l,r}(t), \,r\ge1)$ are just the Poisson probabilities (\ref{plr}).
Hence $1 - p^{l,R}(t) \ge e^{-1}$ if $tp_l \le 1$, and $1 - p^{l,R}(t) \ge
q(r_m) := {\rm Poisson}(1)\{[r_m+1,\infty)\}$ if $tp_l > 1$, implying that
\[
   \min_{l\ge1}(1-p^{l,R}(t)) \ \ge\ c_R \ :=\ \min\{e^{-1},q(m)\} \ >\ 0,
\]
for all~$t$.  It thus follows that $a^T \Sigma_R(t) a \ge c_R a^Ta$ for all $a \in \real^m$.

It is now immediate that, for any $x \in \real^m$, $|\Sigma^{-1/2}_R(t)x| \le c_R^{-1/2}|x|$,
and hence, since $|Y'_{l,r}(t)| \le 1/\sqrt{V_{r}(t)}$ a.s., we have
\[
   |\Sigma^{-1/2}_R(t) (Y'_{l,r_1}(t),\ldots,Y'_{l,r_m}(t))^T|^3
     \Le c_R^{-3/2}\frac{\sjm \{Y'_{l,r_j}(t)\}^2 \sqrt m}{\min_{1\le i\le m} \sqrt{V_{r_i}(t)}};
\]
taking expectations and adding over~$l\ge1$ gives $\b_t \le (m/c_R)^{3/2}/
\min_{1\le i\le m} \sqrt{V_{r_i}(t)}$, proving the first statement of the theorem.
The second follows in view of~\Ref{tv-error}.
\ep

\bigskip
\nin  Thus multivariate normal approximation is always good if the variances of the
(unstandardized) components~$Y_r(t)$ are large.  However, convergence 
typically does not
take place: see a series of  examples in Proposition~\ref{gen-ex} below.

\section{Moments}
\setcounter{equation}{0}
For normal approximation, in view of Theorem~\ref{approximation},
we are particularly interested in conditions under which $V_r(t)\to\infty$.

\par For the moments we have the formulas
\begin{eqnarray}
   \Phi_r(t)&=&\sum_{j=1}^\infty p_{j,r}(t),\label{M}\\
   V_r(t)&=& \sum_{j=1}^\infty p_{j,r}(t) \left(1- p_{j,r}(t)\right)
      \Eq \Phi_r(t) - 2^{-2r} {2r\choose r}\Phi_{2r}(2t),\label{V}\\
   C_{rs}(t) &=& -2^{-r-s}{r+s\choose r} \Phi_{r+s}(2t),\quad r\ne s,\label{cov}
\end{eqnarray}
where, as above,
$p_{j,r} = e^{-tp_j}(tp_j)^r/r!.$

\par 
From   (\ref{M}) and (\ref{V}) we obtain
$$
   \Phi_r(t)\ >\ V_r(t)\ >\ k_r \Phi_r(t),
$$ 
with $k_r > 0$, as is seen from the inequalities
$$ 
    1\ \geq\ 1-{e^{-x}x^r\over r!}\ \geq\  1-{e^{-r}r^r\over r!}\ >\ 0.
$$
for $x\geq 0$.  It follows that 
$$
   V_r(t)\to\infty\ \Longleftrightarrow \ \Phi_r(t)\to\infty;
$$
hence, as long as only the convergence to infinity of~$V_r(t)$ is concerned, 
we can deal with the simpler quantity~$\Phi_r(t)$. This facilitates the
proof of the following theorem, showing how the asymptotic behaviour 
of~$V_r(t)$ for different values of~$r$ is structured.  

\begin{theorem}\label{V-structure}
The asymptotic behaviour of the quantities $V_r(t)$ as $t\to\infty$ follows
one of the following four regimes:
\begin{itemize}
\item
[{\rm 1.}]
  $\lti V_r(t) = \infty$ for all $r\ge1$;
\item
[{\rm 2.}]
  $\limsup_{t\to\infty} V_r(t) = \infty$ for all $r\ge1$, and there exists
  an $r_0\ge1$ such that $\liminf_{t\to\infty} V_r(t) = \infty$ for all $1\le r\le r_0$, 
  and $\liminf_{t\to\infty} V_r(t) < \infty$ for all $r > r_0$;
\item
[{\rm 3.}]
  $\limsup_{t\to\infty} V_r(t) = \infty$ and $\liminf_{t\to\infty} V_r(t) < \infty$ for all $r \ge 1$;
\item
[{\rm 4.}]
  $\sup_t V_r(t) < \infty$ for all $r \ge 1$.
\end{itemize}
\end{theorem}

\nin{\bf Proof.}\ Replacing $V_r$ with~$\Phi_r$ for the argument, the formula~\Ref{M}
yields
\[
   \Phi_r(t) \Eq \sji e^{-tp_j}\frac{(tp_j)^r}{ r!}\,; \qquad 
     \Phi_s(t/2) \Eq \sji e^{-tp_j/2}\frac{(tp_j/2)^s}{s!}\,.
\]
For $s<r$, the ratio of the individual terms is given by
\[
  \frac{e^{-tp_j/2} (tp_j/2)^s/s!}{e^{-tp_j} (tp_j)^r/ r!}
    \ \ge\ \min_{y>0}\{e^{y/2}y^{-(r-s)}\} \frac{r!}{s!2^s} \Eq 
    \left(\frac e{r-s}\right)^{r-s}\,\frac{r!}{s!2^r}.
\]
Hence, for all $s<r$,
\eq\label{doubling-inequality}
   \Phi_s(t/2) \ \ge\ \Phi_r(t)\,\left(\frac e{r-s}\right)^{r-s}\,\frac{r!}{s!2^r}\,.
\en
It now follows that if, for some~$r$, $\lti V_r(t) = \infty$, then $\lti V_s(t) = \infty$
for all $1\le s\le r$ also; and that, if 
$\sup_t V_r(t) < \infty$ for some~$r$, then $\sup_t V_s(t) < \infty$ for all
$s > r$.  Hence, to complete the proof, we just need to show that, if
$\sup_t V_r(t) < \infty$ for some~$r\ge1$, then $\sup_t V_1(t) < \infty$.

For this last part, write $\Phi_r(t) = L_r(t) + R_r(t)$, where
\eq\label{L/R-def}
   L_r(t) \ :=\ \sum_{j\colon tp_j \ge 1} e^{-tp_j}(tp_j)^r/r!\,;\quad
   R_r(t) \ :=\ \sum_{j\colon tp_j < 1} e^{-tp_j}(tp_j)^r/r!\,.
\en
Suppose that $\sup_t\Phi_r(t) = K < \infty$.  Then, for 
every $t>0$,
\eq\label{L-bnd}
  L_1(t) \Le r!\,L_r(t) \Le r!\, \Phi_r(t) \Le K\,r!\,.
\en
It thus remains to bound $R_1(t)$, which in turn can be reduced to finding a bound for
$$
    S(t):=\sum_{j\colon tp_j<1} tp_j.
$$

Let $a_0 \ge a_1\ge\ldots \ge 0$ be any decreasing sequence such that
 $a_j/a_{j+h} \ge 2$ holds for some $h\ge1$ and all $j\ge1$.  
Then $a_{ih+m} \le a_m 2^{-i}$ for every $i\ge 0$ and $0\leq m<h$.
Splitting the $a_j$'s into $h$ subsequences that are dominated by the geometric series,  
we thus have
\[
   \sjo a_j \Le  \sum_{m=0}^h 2 a_m \Le  2a_0 h.
\]
Now if, for some $h\geq 1$, the frequencies~$p_j$ satisfy 
\begin{equation}\label{dense}
   p_j/p_{j+h}\geq 2 \quad \mbox{for all}\quad j\geq 1,
\end{equation}
then applying the above result to the sequence $a_j=tp_{j+\min\{i: tp_i<1\}}$ for any~$t$
yields the bound $R_1(t)<S(t)<2h$, since $a_0 < 1$.

\par On the other hand, if $p_j/p_{j+h}<2$ for some $j$ and $h$, then it follows from 
$p_j\geq p_{j+1}\geq \ldots\geq p_{j+h}>p_j/2$  that  
$$
   L_1(2/p_j)\ >\ \sum_{k=j}^{j+h} e^{-2p_k/p_j} {2p_k\over p_j}\ >\ e^{-2} (h+1).
$$
Thus, for any $h$ such that $e^{-2}(h+1)>Kr!\,$, we see that (\ref{dense}) must hold,
since otherwise (\ref{L-bnd}) would be violated for $t=2/p_j$. Hence it follows that  
$R_1(t)<S(t)<2e^2Kr!\,$,  
and the final part of the lemma is proved.
\ep

\bigskip
\nin In particular, in Theorem~\ref{approximation}, the quantity
$\min_{1\le i\le m} \sqrt{V_{r_i}(t)}$ can thus be replaced in the error 
estimates by $\Phi_{r_m}(2t)$.

\bigskip
We now turn to finding conditions sufficient for distinguishing the
asymptotic behaviour of the~$V_r(t)$.  To do so, introduce the measures 
$$ 
   \nu_r({\rm d}x)=\sum_{j=1}^\infty p_j^r \delta_{p_j}({\rm d}x).
$$
Two special cases are $\nu_0$, a counting measure, and $\nu_1$, the probability 
distribution of a size-biased pick from the~$p_j$'s.
For $r>0$ write (\ref{M}) as
\begin{equation}\label{int-parts}
  \Phi_r(t) \Eq {t^r\over r!} \int_0^\infty e^{-tx} x^r\nu_0({\rm d}x)
  \Eq {t^r\over r!} \int_0^\infty e^{-tx} \nu_r({\rm d}x)
	\Eq {t^{r+1}\over r!} \int_0^\infty e^{-tx} \nu_r[0,x]\,{\rm d}x.
\end{equation}
Comparing with standard gamma integrals, it is then immediate that
\eq\label{gamma}
   \liminf_{x\to0} \frac{\nu_r[0,x]}{x^r} \Le \liminf_{t\to\infty} \Phi_r(t) 
    \Le \limsup_{t\to\infty} \Phi_r(t) \Le \limsup_{x\to0} \frac{\nu_r[0,x]}{x^r}\,.
\en
This, together with Theorem~\ref{V-structure},  enables us to conclude the following
%sufficient 
conditions for the convergence to infinity  of~$\Phi_r(t)$,
and hence equivalently of~$V_r(t)$, expressed in terms of the accessible quantities
$$
    \rho_{j,r}:={1\over p_j^r}\sum_{i=j+1}^\infty p_i^r\,.
$$

\begin{lemma}\label{easy-conds}
\mbox{}

\par
\begin{itemize}
\item[\rm(a)]
$\sup_{t\geq 0} \Phi_s(t)<\infty$ for all $s\ge1$ if and only if, 
for some (and then for all)~$r\ge1$,~ $\sup_{j} \r_{j,r} < \infty$.
\item[\rm(b)]
If, for some~$r\ge1$,
$\lim_{j\to \infty} \r_{j,r}=\infty$,
then~ $\lim_{t\to\infty} \Phi_s(t)=\infty$ for all $1\le s\le r$.
\end{itemize}
\end{lemma}

\nin{\bf Proof.}\ If $p_{j+1}\leq x<p_j$ then 
$$  
   \r_{j,r} \Eq
%   {\nu_r[0,p_{j}]\over p_j^r}-1 \Eq 
   {\nu_r[0,p_{j+1}]\over p_j^r}\Eq
   {\nu_r[0,x]\over p_j^r}\ <\ {\nu_r[0,x]\over x^r} \Le
   {\nu_r[0,p_{j+1}]\over p_{j+1}^r} \Eq 1 + \r_{j+1,r}.
$$
Hence \Ref{gamma} can be replaced by the inequalities
\eq\label{gamma-rho}
  \liminf_{j\to\infty} \r_{j,r} \Le \liminf_{t\to\infty} \Phi_r(t) 
    \Le \limsup_{t\to\infty} \Phi_r(t) \Le 1 + \limsup_{j\to\infty} \r_{j,r}\,.
\en
Part~(b) of the lemma now follows directly from Theorem~\ref{V-structure}.

For part~(a), much as for the last part of the proof of Theorem~\ref{V-structure},
define
\[
   h(j) := \max\{l \ge 0\colon\,p_{j+l}/p_j \ge 1/2\}; \qquad h^* := \sup_j h(j).
\]
Then it is immediate that
\[
      2^{-r}h(j) \Le \r_{j,r} \Le h^*\sli 2^{-(l-1)} \Eq 2h^*,
\]
so that $h^* < \infty$ if and only if $\sup_j \r_{j,r} < \infty$ for some,
and then for all, $r\ge1$.  We now conclude the proof by showing that 
$\sup_{t\geq 0} \Phi_s(t)<\infty$ for all $s\ge1$ if and only if $h^* < \infty$.
Defining $L_r(t)$ and~$R_r(t)$ as in~\Ref{L/R-def}, we observe that, if $h^* < \infty$,
then
\[
   R_r(t) \Le h^* \sli 2^{-r(l-1)} \Le 2h^* \quad\mbox{and}\quad
   L_r(t) \Le h^* \sli e^{-2^{l-1}}\,\frac{2^{lr}}{r!}\,,
\]
so that $\Phi_r(t) = L_r(t) + R_r(t) < \infty$ for all~$r\ge1$.  On the other hand,
\[
   L_r(1/p_{j+h(j)}) \ \ge\ e^{-2}h(j)/r!\,,
\]
implying that, if $h^* = \infty$, then $\limsup_{t\to\infty}\Phi_r(t) = \infty$
for all $r\ge1$.
\ep

\bigskip
The familiar ratio test yields simpler sufficient conditions. Thus $\sup_t \Phi_r(t) < \infty$ 
for all $r\ge1$ if 
$$
    \limsup_{j\to\infty} p_{j+1}/p_j\ <\ 1,
$$ 
while $\lti \Phi_r(t) = \infty$ for all $r\ge1$ if 
$$
    \lim_{j\to\infty} p_{j+1}/p_j \Eq 1.
$$
% in this latter case, the successive $p_j$'s are almost equal, whence $\rho_{j,r}\to\infty$. 
For instance, for $p_j=cq^j$, the geometric distribution with $0<q<1$, 
we have $p_{j+1}/p_j=q$; hence
$\sup_t \Phi_r(t)<\infty$ for all~$r$, and normal approximation is not adequate
for any~$r$.  This illustrates possibility~4 in Theorem~\ref{V-structure}.
For the Poisson distribution $p_j=c\lambda^j/j!$\,, we even have $p_{j+1}/p_j\to 0$, 
and so normal approximation is no good here, either.

\par Continuing this line, we obtain a further set of conditions.

\begin{lemma}\label{harder-conds}
\begin{itemize}
\item[\rm(a)] Suppose for some $0<\lambda<1$
\begin{equation}\label{coco1}
\liminf_{j\to\infty} {p_{j+h}\over p_j}>\lambda
\end{equation}
for every $h\geq 1$. 
Then $\Phi_r(t)\to\infty$ as $t\to\infty$ for all $r\geq 1$.
\item[\rm(b)] The condition $\limsup_{t\to\infty} \Phi_r(t)<\infty$ holds for some (hence for all) $r\geq 1$ if and only if
there exists $h\geq 1$ such that 
\begin{equation}\label{coco2}
\limsup_{j\to\infty} {p_{j+h}\over p_j}\leq {1\over 2}.
\end{equation}
\end{itemize}
\end{lemma}
\nin{\bf Proof.}\ For part (a), assume that $\nu_0(\lambda x,x)=\#\{j: \lambda x<p_j<x\}\to\infty$ as $x\to0$.
Then also
$$\Phi_r(1/x)\geq \sum_{\{j: \lambda x<p_j<x\}} e^{-p_j/x} (p_j/x)^r/r!\geq \nu_0(\lambda x,x)\,\min_{\{y:\lambda<y<1\}}[e^{-y}y^r/r!]
\to\infty.$$ 
As $x$ decreases, the piecewise-constant function $\nu_0(\lambda x,x)$ may have downward jumps only at the values
$x\in \{p_j\}$, hence the assumption is equivalent to $\nu_0(\lambda p_j,p_j)\to\infty$ (as $j\to\infty$), which 
in turn is readily translated into (\ref{coco1}).

\par For part (b), the same estimate with any $0<\lambda<1/2$ shows that the condition (\ref{coco2}) 
is necessary.
In the other direction, suppose that $p_{j+h}/p_j<3/4$ for all $j\geq J$. Split $(p_j, j\geq J)$ into 
$h$ subsequences 
$(p_{J+s+ih}, i\geq 0)$,  with $0\leq s\leq h-1$. Each of the subsequences has the property that 
the ratio of any two consecutive elements is at most $3/4$.
Hence, as above, the sum of the terms $e^{-p_jt}(tp_j)^r/r!$  along a subsequence
yields a  uniformly bounded contribution to $\Phi_r$.
\ep

\bigskip
\par Examples of irregular behaviour of moments may be constructed by breaking  
the sequence $(p_j,\,j\ge1)$ into finite blocks of sizes
$m_1,m_2,\ldots$, and setting the $p_j$'s within the $i$'th block all equal to some $q_i$.
We use the notation $V(t) := \var\bigl(\sum_{r\ge1}Y_r(t)\bigr)$ to denote the
variance of the number of occupied boxes.

\bigskip
\nin {\bf Example~1.} \cite[p. 384]{Karlin}. Take $m_i=i$ and $q_i=c 2^{-2^i}$, 
with $c$ a normalizing factor\footnote{In fact, the Poisson sampling model
makes sense for arbitrary $p_j$'s, and the enumeration of small counts 
makes sense if $\sum_j p_j<\infty$.} 
to achieve $\sum_j p_j=1$.
Then both $V(t)$ and $\Phi_1(t)$ oscillate between $0$ and $\infty$, 
approaching the extremes arbitrarily closely.  This illustrates possibility~3 in
Theorem~\ref{V-structure}.

\bigskip
\nin {\bf Example~2.} As in \cite[Example 4.4]{BGY}, take $m_i=2^{2^i}$, $q_i=c2^{-2^{i+1}}$. 
Then $\Phi_1(t)\to\infty$, but $\Phi_2(t)$ oscillates between $0$ and~$\infty$
as $t$ varies; thus $Y_1(t)$ is asymptotically normal, but $Y_2(t)$ is not, and
% This shows that~(b) cannot be understood as equivalence, at least for $r=1$. 
the ratios $p_{j+1}/p_j$ have accumulation points at $0$ and $1$.
This illustrates possibility~2 in
Theorem~\ref{V-structure}.

\bigskip
\nin We now extend this example, showing among other things that one can have 
any value for~$r_0$ in behaviour~2 in Theorem~\ref{V-structure}.

\begin{proposition}\label{gen-ex}
Fix $0 < \b < 1$ and $\a > 0$, and
take the blocks construction with $m_i = \lfloor 2^{(1-\b)^{-i}}\rfloor$, $q_i = cm_i^{-(1+\a)}$,
where~$c$ is the appropriate normalizing constant.  Then we have
\begin{itemize}
\item[\rm(i)]
  $\limsup_{t\to\infty} V_r(t) \Eq \infty$ for all $r \ge 1$;
\item[\rm(ii)]
   $\lti V_r(t) \Eq \infty$ if and only if $r\b(1+\a) \le 1$;
\item[\rm(iii)]
  $\lim_{j\to\infty}\r_{r,j} = \infty$ if and only if $r\b(1+\a) < 1$;
\item[\rm(iv)]
  The quantities $\Sigma_{rs}(t)$ do not converge for any $r\ne s$.
\end{itemize}
\end{proposition}

\nin{\bf Proof.}\ Once again, we work with $\Phi_r$ instead of~$V_r$, now writing
\eq\label{f-sum}
    \Phi_r(t) = \sii m_i e^{-tq_i}(tq_i)^r/r!\,.
\en
For part~(i), it is enough to consider the subsequence $t_l := 1 / q_l$, $l\ge1$.

For part~(ii), split~$\real_+$ into intervals $J_l := [q_l^{-1},q_{l+1}^{-1})$, $l\ge1$;
we show that $\lim_{l\to\infty}\inf_{t\in J_l} \Phi_r(t) = \infty$ if $r\b(1+\a) \le 1$,
and exhibit a subsequence $(t'_l,\,l\ge1)$ with $t'_l \in J_l$ such that $\lim_{l\to\infty} 
\Phi_r(t'_l)=0$ if $r\b(1+\a) > 1$. Indeed, for $t\in J_l$, taking just the term with $i =l+1$
in~\Ref{f-sum}, we obtain
\[
    m_{l+1} \exp\{-\f q_{l+1}/q_l\} (\f q_{l+1}/q_l)^r/r!
   \ \asymp\ m_{l+1} \f^r \left(\frac{m_{l+1}^{(1-\b)(1+\a)}}{m_{l+1}^{(1+\a)}}\right)^r
     \Eq \f^r m_{l+1}^{1-r\b(1+\a)},
\]
where we write $t = \f/q_l$ with $1 \le \f \le q_l/q_{l+1} \sim m_{l+1}^{\b(1+\a)}$, and
use the fact that $\f q_{l+1}/q_l \le 1$ in this range.  For $r\b(1+\a) < 1$, it follows
that $\inf_{t\in J_l} \Phi_r(t) \asymp m_{l+1}^{1-r\b(1+\a)} \to \infty$ as $l\to\infty$.

For $r\b(1+\a) = 1$, take also the term with $i = l$ in~\Ref{f-sum},
giving a combined contribution of at least
\[
   m_l e^{-\f}\frac{\f^r}{r!} + K\f^r,
\]
for some $K > 0$. It is easily checked that the minimum value of this sum for $\phi>1$ goes to $\infty$ with $l$,
hence, once again, $\lim_{l\to\infty}\inf_{t\in J_l} \Phi_r(t) = \infty$.

For $r\b(1+\a) > 1$, these two terms contribute an amount of order
\eq\label{rba-big}
   \f^r\{m_{l+1}^{(1-\b)} e^{-\f} + m_{l+1}^{1-r\b(1+\a)}\},
\en
to~\Ref{f-sum}, which is small as $l\to\infty$, for example, for $\f = 2\log m_{l+1}$. 
The sum of the terms in~\Ref{f-sum} for $i \ge l+2$ is of order
\[
   \sum_{i \ge l+2} m_i \left(\frac{\f q_i}{q_{l}}\right)^r
     \ \sim\  \f^r \sum_{i \ge l+2} m_i\{m_i^{(1-\b)^{i-l}-1}\}^{r(1+\a)}
     \Eq \f^r O(m_{l+1}^{1-r\b(1+\a)-\h}),
\]
where $\h > 0$, and hence asymptotically smaller than the second element of~\Ref{rba-big}.
The sum of the terms in~\Ref{f-sum} for $i \le l-1$ is of order at most
\[
   \left\{\sum_{i=1}^{l-1} m_i\right\} \exp\{-\f q_{l-1}/q_l\} 
    \left(\frac{\f q_{l-1}}{q_l}\right)^r,
\]
largest for $\f=1$ for all~$l$ large enough, when it is of order
\[
    m_{l-1}^{1 + r\b(1+\a)/(1-\b)} \exp\{-m_{l-1}^{\b(1+\a)/(1-\b)}\},
\]
asymptotically small as $l\to\infty$.  Hence, for $t'_l = 2q_l^{-1}\log m_{l+1}$,
it follows that $\lim_{l\to\infty} \Phi_r(t'_l) = 0$, and therefore that $\Phi_r(t)$
does not converge to infinity as $t\to\infty$.

For part~(iii), writing $M_i := \sum_{l=1}^i m_l$, we have 
\[
   \r_{r,j} \ \ge\ q_i^{-r}\sum_{l\ge i+1} m_l q_l^r\quad 
      \mbox{whenever}\quad M_{i-1} < j \le M_i,
\]
with equality for $j=M_i$.  Now
\[
    \sum_{l\ge i+1} m_l q_l^r \ \asymp\ m_{i+1}^{1-r(1+\a)},
\]
and
\[
    q_i^{-r} \Eq m_i^{r(1+\a)} \ \sim\ m_{i+1}^{r(1-\b)(1+\a)}.
\]
Hence $\r_{r,M_i} \asymp m_{i+1}^{1-r\b(1+\a)}$ is bounded for $r\b(1+\a) \ge 1$,
and $\r_{r,j} \to \infty$ as $j\to\infty$ if $r\b(1+\a) < 1$.

For part~(iv), we 
%use the same intervals as were introduced for part~(ii), and 
note
that, for $t = \f / q_l$, 
%with~$\f$ restricted to any fixed interval $[1,F]$,
the quantity
\[
   \Sigma_{rs}(t) = -2^{-r-s} {r+s\choose r}\ \frac{\Phi_{r+s}(2t)}{\sqrt{V_r(t)V_s(t)}},\quad r\ne s,
\]
behaves asymptotically, as~$l$ becomes large, in the same way as for the Poisson occupancy 
scheme with a single block
of $m_l$ boxes with equal frequencies $q_l$. Computing the limit,
\[
\lim_{l\to\infty} \Sigma_{rs}(\phi/q_l)=\,\, -{1\over \sqrt{r!s!}} {e^{-\phi}\phi^{(r+s)/2}\over  
\sqrt{\{1-e^{-\f}/r!2^r\}\{1-e^{-\f}/s!2^s\}}}\,,
\]
where $m_l$ cancels because of the additivity of the moments. As $\phi$ varies, this 
limit value varies too, and hence, for $r\ne s$,
the quantities $\Sigma_{rs}(t)$ do not converge as $t\to\infty$.
\ep

\bigskip
It follows from parts (ii) and~(iii) of Proposition~\ref{gen-ex}
that the implication in part~(b) of
Lemma~\ref{easy-conds} cannot be reversed, and from part~(iv) that the correlations
between different components of $Y(t)$ need not converge, even when their variances
tend to infinity.  Hence the approximation in Theorem~\ref{approximation} does
not necessarily imply convergence.
Yet another kind of pathology appears when $Y_1(t)$ is asymptotically independent of $(Y_r(t), r>1)$, as in the following example.

\bigskip
\nin {\bf Example~3.} 
Suppose that the
frequencies in the block construction satisfy  
% $m_i\to\infty$, ${q_{i+1}/ q_i}\to 0$ and ${q_{i+1}m_{i+1}/(q_im_i})\to 1$ as $i\to\infty$. 
% A possible choice is 
$q_i=1/i!,~m_i=(i-2)!$ (with $i\geq 2$). Since
$q_i^{-r}\sum_{k=i+1}^\infty m_kq_k^r \to \infty$ for each~$r$,
we have $\lim_{j\to\infty} \r_{j,r} = \infty$, and hence
all the variances $V_r(t)$ go to $\infty$ by Lemma~\ref{easy-conds}\,(b). 
On the other hand, $m_iq_i\Big/\sum_{k=i+1}^\infty m_kq_k \to 0$, and it follows that
$$
   {\Phi_{1+s}(2t)\over\Phi_1(t)}\ =\ 
     {2^{s+1}\sum_i m_iq_ie^{-tq_i} \{e^{-tq_i} t^sq_i^s\}\over (s+1)!\,\sum_i m_iq_i e^{-tq_i}}
        \ \to\ 0
$$
as $t\to\infty$.  Since $\Phi_{1+s}(2t)/\Phi_s(t)$ is bounded above by~\Ref{doubling-inequality},
we conclude that $\Sigma_{1,s}(t)\to 0$ for $s\ge2$. It follows that every pair $(Y_1'(t),Y_s'(t))$, $s\geq 2$,
converges in distribution to the standard bivariate normal distribution with independent components.
Because the variances go to $\infty$, Theorem \ref{approximation} guarantees increasing quality 
of the normal approximation for any finite collection of components $Y_{r_i}'(t)$. 
However, the full vector $(Y_r', r=1,2,\ldots)$  does not converge: see more on this example in Sections 5 and 6.

\bigskip

\par Part (ii) of  Proposition~\ref{gen-ex} also demonstrates that $\liminf_{j\to\infty} p_{j+1}/p_j=0$ 
does not exclude 
that $\Phi_r(t)\to\infty$, hence the condition (\ref{coco1}) in Lemma \ref{harder-conds} is not necessary.
Finally, by \cite[Eqn. 3.1]{BGY}, we have
$$
   {1\over 2} \Phi_1(2t)\ <\ V(t)\ <\ \Phi_1(t),
$$
meaning that $\Phi_1(t)$ is always of the same order as the variance of the number 
of occupied boxes $V(t)$.  The examples above show that this need not
be the case for $\Phi_r(t)$, when $r\ge2$.

\section{Regular variation}\label{RV}
\setcounter{equation}{0}
We now henceforth assume that $\Phi_r(t)\to\infty$ for all $r\geq 1$. 
The CLT for each component of~$Y_{t}$ then holds, as observed above, and normal approximation
becomes progressively more accurate for the joint distribution of any finite
collection of components. A joint normal limit for any collection of
the standardized components also holds, provided that  the corresponding 
covariances converge. From (\ref{cov}) we have
\begin{equation}\label{cov2}
    {\rm Cov\,}(Y_{r}'(t), Y_{s}'(t)) \Eq \Sigma_{rs}(t) \Eq c(r,s)
      {\Phi_{r+s}(2t)\over \sqrt{V_{r}(t)V_{s}(t)}},~~~~~r\neq s.
\end{equation}
The RHS  converges to a nonzero limit for each pair $r,s$ if, for each~$r$,  $\Phi_r\approx f
\in R_\alpha$, where $R_\alpha$ denotes the class of functions regularly varying at $\infty$
 with index $\alpha$, and where, here and subsequently, we write
$a\approx b$ if $a(t)/b(t)\to c$ as $t\to\infty$ with $0<c<\infty$.
If $\Phi_r\in R_\alpha$, then the index belongs to the range $0\leq\alpha\leq 1$, 
because $\Phi_r(t)$ cannot converge to $0$, and because $\Phi_r(t)/t\to 0$.

\par The results in the next section show that, if the covariances converge for a 
sufficiently large set of pairs $r,s$, then this is in fact the only possibility.
More formally,
we say that then regular variation holds in the occupancy problem,  meaning that,
for some $0\leq\alpha\leq 1$ and some {\em rate function\/} $f \in R_\alpha$, 
\begin{equation}\label{RVnew}
   \Phi_r\ \approx\ f  ~~~{\rm for~all~~} r\geq 2\,.
\end{equation}
This setting of regular variation extends the original approach by Karlin \cite{Karlin} in 
the  special case $\alpha=0$,
and, moreover, it covers all possible limiting covariance structures (Theorem \ref{forceRV}).

\par Observe that the functions $t^{-r}\Phi_r$ satisfy 
\begin{equation}\label{one-pr}
{{\rm d}^r\over{\rm d}t^r}\left\{t^{-1}\Phi_1(t)\right\}=(-1)^rr!\left\{t^{-r}\Phi_r(t)\right\},
\end{equation}
thus, in particular, they are completely monotone. 
This taken together with the standard properties  of regularly varying functions  \cite{BGT} 
implies that, if $\Phi_r\in R_\alpha$ for some $0\leq \alpha<1$ and
$r\geq 1$, then the same is true for all $r\geq 1$, and we can choose the rate function $f=\Phi_1$. 
The case $\alpha=1$ is special. If $\Phi_r\in R_1$ for some $r\ge2$, then all $\Phi_r$ for $r\ge2$ are 
of the same order of growth and $\Phi_1\in R_1$, but $\Phi_1\gg \Phi_2$
(this motivates the choice $r\ge2$ in~\Ref{RVnew}).

\par A necessary condition for \Ref{RVnew} is  $\lim_{j\to\infty} p_{j+1}/p_j=1$,
as follows from the next lemma.

\begin{lemma}\label{noRV}
If~ $\liminf_{j\to\infty} p_{j+1}/p_j<1$ then $\Phi_r$ is not regularly varying  for $r\geq 2$,
and $\Phi_1$ is not regularly varying with index $\alpha<1$.
\end{lemma}
\nin{\bf Proof.}\
We have 
$$
   t^{-2}\Phi_2(t)=\sum_{j=1}^\infty e^{-tp_j}p_j^2=\int_0^1 e^{-tx}\nu_2({\rm d}x)
$$
with $\nu_2[0,x]:=\sum_{j=1}^\infty p_j^2 \,{\bf 1}[p_j\leq x]$. 
Suppose $t^{-2}\Phi_t \in R_{-\beta}$, then $1\leq \beta\leq 2$ and, by Karamata's Tauberian 
theorem, also  $\nu_2[0,t^{-1}]\in R_{-\beta}$.
Because $\beta\neq 0$, the latter implies that $\nu_2[at^{-1},bt^{-1}]\in R_{-\beta}$, i.e. that 
\begin{equation}\label{disp}
\nu_2[at^{-1},bt^{-1}]\sim (b^\beta-a^\beta)\ell(t)t^{-\beta},\quad t\to\infty
\end{equation}
for any positive $a<b$. However, the assumption of the lemma allows to choose $a<b<1$ such 
that $\nu_2[ap_j,bp_j]=0$ for infinitely many $j=j_k$,
so (\ref{disp}) fails for $t=1/p_{j_k}\to\infty$. The contradiction shows that $t^{-2}\Phi_2(t)$ 
cannot be regularly varying.
The assertions regarding $r\neq 2$ can be derived in the same way.
\ep

\bigskip
The example below shows that $\Phi_r$ may be regularly varying for $r=1$ alone.

\medskip
\noindent
{\bf Example 3} (continued). 
%In this example $\Phi_1$ is of regular variation only for $r=1$.
Let $g(t)=\nu_1[0,t^{-1}]=\sum_{j=1}^\infty p_j\,{\bf 1}[p_j\leq t^{-1}]$.
We have the general estimates
$$
    t^{-1}\Phi_1(t)\ \ge\ e^{-1} g(t)
$$
and, for $a > 1$ and any $\epsilon > 0$,     
\begin{eqnarray*}
     \lefteqn{t^{-1}\Phi_1(t) - (at)^{-1}\Phi_1(at)}
      \\
  &&\Le \epsilon g(at/\epsilon) + 
    \{g(t/\log\{1/\epsilon g(t)\})-g(at/ \epsilon)\}
               + \sum_{j=1}^\infty p_j e^{-tp_j}  \,  {\bf 1}[ p_j > t^{-1}\log\{1/\epsilon g(t)\}]
     \\
  &&\Le 2\epsilon g(t) + \{g(t/\log\{1/\epsilon g(t)\}) - g(at/\epsilon)\}.
\end{eqnarray*}
Applying these to the block construction with $q_i=1/i!$ and $m_i=(i-2)!$, we observe that 
$g(t) \asymp I(t)^{-1}$ and that $g(t/\log\{1/\epsilon g(t)\}) - g(at/\epsilon)$ involves 
at most two $q_i$, each 
of the corresponding terms being of the order of $I(t)^{-2}$,
where  $I(t):=\min\{i:~i!\geq t\}$.
  It follows that $t^{-1}\Phi_1(t)\in R_0$, 
whence $\Phi_1\in R_1$ and $\Phi_1\gg \Phi_r$ for $r\ge2$. However, $q_{i+1}/q_i\to 0$, therefore
Lemma \ref{noRV} implies that $\Phi_r\notin R_1$ for $r\ge2$.

\vskip0.5cm

\par The {\it proper} case of regular variation  with index $0<\alpha<1$  
can be characterized by Karlin's condition \cite[Equation 5]{Karlin}
\begin{equation}\label{regvar}
    \nu_0[x,1]:=\#\{j: p_j\geq x\} \sim\ \ell(1/x)x^{-\alpha}, ~~~~~~~x\downarrow 0,
\end{equation}
where and henceforth
the symbol $\ell$ stands for a function of slow variation at $\infty$. Other equivalent conditions are 
(see \cite{GHP})
\begin{eqnarray*}
\Phi(t)&:=&\int_0^1(1-e^{-tx})\nu_0({\rm d}x)\sim \Gamma(1-\alpha)t^\alpha\ell(t),\\
\nu_r[0,x]& \sim & {\alpha\over r-\alpha}\, x^{r-\alpha}\ell (1/x) {\rm~~~ for~ some~}r\geq 1,\\
\Phi_r(t) &\sim &  {\alpha\Gamma(r-\alpha)\over r!}\,t^\alpha\ell(t){\rm~~~ for~ some~}r\geq 1,\\
p_j &\sim & \ell^*(j) j^{-1/\alpha},
\end{eqnarray*}
where $\ell^*(y)=1/\{\ell^{1/\alpha}(y^{1/\alpha})\}^{\#}$, and $\#$ denotes the de Bruijn conjugate
of a slowly varying function \cite{BGT}.
Note that $V_r(t)$ then has the same order of growth, in view of~\Ref{V}, yielding
behaviour as in possibility~1 of Theorem~\ref{V-structure}. 
The joint CLT for
$${Y_r(t)-\Phi_r(t)\over \sqrt{t^\alpha\ell(t)}},~~~~~r=1,2,\ldots$$
 in~$\real^\infty$ holds with the limiting covariance matrix $S$ computable from (\ref{cov}) as
\begin{eqnarray*}
S_{rs}= -{\alpha\Gamma(r+s-\alpha)\over r!s! 2^{r+s-\alpha}}\,,~~~~~r\neq s\\
S_{rr} ={\alpha\over r!}\left(\Gamma(r-\alpha)-{\Gamma(2r-\alpha)\over r!2^{2r-\alpha}}\right),
\end{eqnarray*}
in accord with Karlin \cite[Theorem 5]{Karlin}.

\par If (\ref{regvar}) holds with $\alpha=1$
then $\ell(t)$ must approach $0$ as $t\to\infty$
sufficiently fast to have $\sum p_j<\infty$\footnote{
One example is $p_j=c/j\{\log (j+1)\}^{\beta+1}$, $\beta>0$, in which case $\ell(t)\sim 1/c(\log t)^{\beta+1}$.}.
In this situation
 we have 
$\Phi_r(t)\sim  (r^2-r)^{-1}\ell(t)t$ for $r>1$ but
$\Phi_1(t)\sim \ell_1(t)t$ with some $\ell_1\gg\ell$. In fact, $X_{n,1}\sim K_n$ as $n\to\infty$ 
almost surely.
Because the scaling of $Y_1(t)$ is faster than that for other $Y_r(t)$'s, 
it follows from \Ref{V} and~\Ref{cov} that $\Sigma_{1r}(t) \to 0$ for all $r\ge2$,
so that the CLT 
holds with $Y_1'(t)$ asymptotically independent of $(Y_r'(t), r\ge2)$.
The limiting covariance matrix of $\{(Y_r-\Phi_r)/(t\ell(t)),~ r\geq 2\}$ is obtained by 
setting $\alpha=1$ in the above formulas for $S$.
Our multivariate result extends in this case the 
 marginal convergence that was stated in 
 \cite[Thm $5'$]{Karlin}\footnote{
Mikhailov \cite{Mikhailov} indicated yet other situation where the $X_{n,r}$'s 
for $r>1$ all behave similarly, but their behaviour is distinct from that of  
 $X_{n,1}$.}.

\par   Karlin's condition  (\ref{regvar}) with  $\alpha=0$ 
is too weak to control the $\Phi_r(t)$'s. However,
a slightly stronger condition 
\begin{equation}\label{regvar1}
    \nu_1[0,x]:=\sum_{\{j:p_j\leq x\}} p_j\ \sim\  x \ell_1(1/x),
\end{equation}
is equivalent to $\Phi_r\in R_0$ for any (and hence for all) $r\geq 1$.
To illustrate the difference, note that in the geometric case, with 
$p_j=(1-q)q^{j-1}$,  $0<q<1$, we have $\ell(1/x)\sim \log_q(1/x)$, whereas 
$\nu_1[0,x]=q^{\lceil \log_q (x/(1-q)) \rceil}$ is not regularly varying, since 
$\nu_1[0,x]/x$ jumps infinitely often from $(1-q)^{-1}$ to $q(1-q)^{-1}$ as $x\to0$.
The geometric case can be contrasted to the one with frequencies
$p_j= ce^{-j^\beta}$ ($0<\beta<1$), for which we have $\ell(1/x)\sim c|\log x|^{1\over\beta}$ and 
$\nu_1 [0,x]/x \sim c |\log x|^{{1\over\beta}-1}$.

 By \cite[Prop. 15]{GHP}, the general connection between $\ell_1$ in (\ref{regvar1}) 
and $\ell$ in (\ref{regvar}) is
$$
  \ell(1/x)\Eq\int_{x}^1 u^{-1}\ell_1(1/u)\,{\rm d}u\,,~~~~~~~0<x<1.
$$
Adopting (\ref{regvar1}) we have
$\nu_r[0,x]\sim r^{-1}x^r\ell_1(1/x)$, $r\ge1$, and the situation is then very similar to 
that in the proper case: we have $\Phi_r(t)\sim r^{-1}\ell_1(t)$ and
$\{(Y_r(t)-\Phi_r(t))/\sqrt{\ell_1(t)},\ r\ge1\}$,
converges in law to a multivariate Gaussian limit with covariance matrix~$S$ given by
$$
  S_{rr}\Eq\left({1\over r}-{1\over r \,2^{2r+1}}{2r\choose r}\right),
       ~~~~S_{rs}\Eq-{1\over (r+s) 2^{r+s}}{r+s\choose r},~~~~r\neq s.
$$
This applies, for instance, to the frequencies $p_j\sim ce^{-j^\beta}$ ($0<\beta<1$). 
This case of slow variation seems not to have been considered before.

\section{Convergence of the covariances}
\setcounter{equation}{0}

We will show in this section that regular variation is essential for the multivariate convergence
of the whole standardized vector of counts, so that all possible 
limit covariance structures are those characterized in the previous section.
Our starting point is the  following  lemma,
which asserts that the regular variation is forced by the  
convergence of the ratios of $\Phi_r$'s.

\begin{lemma} Suppose for some $r\geq 1$
\begin{equation}\label{ratio1}
    \lti \Phi_{r+1}(t)/\Phi_r(t) \Eq c.
\end{equation}
Then $(r-1)/(r+1)\leq c\leq r/(r+1)$ and $\Phi_r\in R_\alpha$ with $\alpha:=r-c(r+1)$. 
Moreover, we then always have 
\begin{equation}\label{ratio2}
    \lim_{t\to\infty}{\Phi_s(t)\over \Phi_r(t)}= {r! \,\Gamma(s-\alpha) \over s!\,\Gamma(r-\alpha)}
\end{equation}
and $\Phi_s\in R_\alpha$ for all $s\geq 1$, unless
$\alpha=1$. If {\rm (\ref{ratio1})} holds with $r>1$ and $c=(r-1)/(r+1)$, then $\Phi_s\in R_1$ for $s\ge2$,
and {\rm (\ref{ratio2})} is still true (in particular, $\Phi_1\gg \Phi_2$).

\end{lemma}
\nin{\bf Proof.}\ A monotone density result which dates back to von Mises and 
Lamperti \cite{Lamperti} says that 
the convergence $tg'(t)/g(t)\to \beta$ implies $g\in R_\beta$
(this holds for arbitrary $\beta$, including $\pm\infty)$.
This result applied to $g(t)=t^{-r}\Phi_r(t)$ yields the regular variation $\Phi_r\in R_\alpha\,$, 
with some  $0\leq\alpha\leq 1$.
 The rest follows from (\ref{one-pr}), monotonicity and the general behaviour 
of the regularly varying functions under integration and differentiation \cite{BGT}.
\ep 

\medskip
\par To apply the lemma, we need  to pass from the convergence of covariances (\ref{cov2}) 
to the convergence of a ratio as in (\ref{ratio1}).
To this end, it is useful to exclude zero limits.

\begin{lemma}\label{nozeroes} 
If $\limsup_t\Phi_s(t) = \infty$ for any~$s\ge1$,
then no correlation $\Sigma_{r,r'}(t)$ with $2 \le r < r'$ can converge to zero. 
\end{lemma}

\nin{\bf Proof.}\
\noindent(i)  Let $m_j := \#\{l:\,2^{-(j+1)} < p_l \le 2^{-j} \}$.  Then, if $m^* := \sup_j m_j < \infty$,
it follows that, for $2^j  \le t < 2^{j+1}$,
\begin{eqnarray*}
   s! \Phi_s(t)   & = &  \sum_{k\ge0}~~\sum_{\{l: 2^{-(k+1)} < p_l  \le 2^{-k}\}} (tp_l)^s e^{-tp_l} \\
       & \le &  \sum_{k\ge0} m_k 2^{(j+1-k)s} \exp\{-2^{j-k-1}\} \\
       &  \le & m^* \Bigl( \sum_{k\ge j+1} 2^{(j+1-k)} + 2^s \sum_{k=0}^j 2^{s(j-k)}\exp\{-2^{j-k-1}\} \Bigr)\\
       & \le &  m^*\Bigl( 2 + 2^s\sum_{l\ge0} 2^{ls}\exp\{-2^{l-1}\}\Bigr) \ =\ m^* c_s \ <\ \infty,
\end{eqnarray*}
uniformly  in~$j$, which contradicts $\limsup_t\Phi_s(t) = \infty$.      Hence $\sup_j m_j = \infty$.
\medskip

\noindent(ii)  Given any $j_0$, there exists some $j \ge j_0$ such that
\begin{equation}\label{requir}
   m_k \le m_j,\quad 0\le k\le j; \qquad m_k \le 3^{k-j} m_j,\quad k\ge j. 
\end{equation}
To see this, first take $j_1 \ge j_0$ such that $m_{j_1} = \max\{m_k,\,0 \le k\le j_1\}$, as
can always be done, since $\sup_j m_j = \infty$.  Then let $j_2 := \max\{k \ge j_1:
m_k \ge 3^{k-j_1} m_{j_1}\}$; this is finite, since $1 \ge \sum_{l\ge1}p_l \ge m_j 2^{-(j+1)}$
for each $j\ge0$.  Finally, take $j_3 = \arg\max_{j_1 \le j \le j_2} m_j$; then $j_3$
satisfies the requirements of (\ref{requir}).
\medskip

\noindent(iii)  Now suppose that $j$ satisfies (\ref{requir}).  Then, much as in
part~(i), for any $r\ge2$,
\begin{eqnarray*}
   r! \Phi_r(2^j) &\Le &  \sum_{k\ge0} m_k 2^{(j-k)r} \exp\{-2^{j-k-1}\} \\
       & \Le & \Bigl( \sum_{k\ge j+1} m_j 3^{k-j}2^{r(j-k)} + m_j 
                               \sum_{k=0}^j 2^{r(j-k)}\exp\{-2^{j-k-1}\} \Bigr)\\
        & \Le & m_j\Bigl( 3 + \sum_{l\ge0} 2^{lr}\exp\{-2^{l-1}\}\Bigr) \ =\ c'_r m_j,
\end{eqnarray*}
with $c'_r < \infty$,
whereas also, just from the indices~$l$ with $2^{-(j+1)} < p_l \le 2^{-j}$, we have
$$
    r! \Phi_r(2^{j+1}) \ \ge\   m_j e^{-2}.
$$
This implies that 
$$
    \Phi_{r+r'}(2t) / \sqrt{\Phi_r(t)\Phi_{r'}(t)}\  \ge\  {e^{-2}/\{r+r'\}! \over \sqrt{c'_rc'_{r'}/r!r'!}}
      \  >\  0
$$
for $t = 2^j$, whenever~$j$ satisfies the requirements of~\Ref{requir}, and there are
infinitely many such.  Hence the correlations $\Sigma_{r,r'}(t)$ with $r' > r \ge 2$ cannot
converge to zero. 
\ep
\medskip

\par Note that the correlations $\Sigma_{1,s}(t), s>1,$ converge to zero in the case of regular variation
with index $\alpha=1$. Example 3 illustrates that $\Sigma_{1,s}(t)$ may also converge to zero 
when regular variation in the sense of~\Ref{RVnew} does not hold.

\bigskip
\begin{lemma}\label{sqrtlemma}  If~$g$ is continuous and positive, and $g(2t)/\sqrt{g(t)} \to k$ 
as $t\to\infty$, with $0 < k < \infty$, then $g(t) \to k^2$.
\end{lemma}
\medskip
\nin{\bf Proof.}\  Given $\e > 0$, let $t_\e$ be such that $g(2t) \le k\sqrt{(1+\e)g(t)}$
for all $t \ge t_\e$.  Let $K_\e := \sup_{t \in J_\e} g(t)$, where $J_\e := [t_\e , 2t_\e]$.
Then, for all $t \in J_\e$ and all $n\ge0$, we have
$$
     g(2^nt) \ \le\ \{k^2(1+\e)\}^{1 - 2^{-n}} \{g(t)\}^{2^{-n}} \ \le\ k^2(1+\e) K_\e^{2^{-n}}.
$$
Thus $\limsup_t g(t) \le k^2$.  A similar argument shows that $\liminf_t g(t) \ge k^2$,
proving the lemma.      
\ep
\bigskip 

\begin{theorem}\label{forceRV}  Suppose the correlations $\Sigma_{r,s}(t)$ converge, as 
$t\to\infty$, for $r,s$ satisfying $2\leq r<s$ and $r+s\leq 12$.
Then  the following is true: 
\begin{itemize}
\item[{\rm (i)}]
 {\rm (\ref{RVnew})} holds with some $0\leq\alpha\leq 1$,
\item[{\rm (ii)}]
 the  correlations $\Sigma_{r,s}(t)$  converge for all $r,s$,
 \item[{\rm (iii)}] 
$(Y_r'(t), r=1,2,\ldots)$ converges weakly to one of the multivariate 
normal laws described in {\rm Section \ref{RV}},
\item[{\rm (iv)}] the same multivariate normal limit holds for the normalized and centred $X_n$.
\end{itemize}
\end{theorem} 
\medskip
\nin{\bf Proof.}\ For short, write $V_j=V_j(t)$, $f_j=\Phi_j(t)$ and $F_j=\Phi_j(2t)$.

\par By Lemma \ref{nozeroes}, the $\Sigma_{r,s}(t)$ converge to nonzero limits, whence, for $r,s$ in 
the required range,
$$ 
    {F_{r+s}\over \sqrt{V_rV_s}}\ \approx\ {F_{r+s}\over \sqrt{V_{r+1}V_{s-1}}}
$$
and hence $V_rV_s\approx V_{r+1}V_{s-1}$. From this, 
$V_5\approx V_3V_4/V_2,~V_6\approx V_3V_5/V_2\approx V_3^2V_4/V_2^2$,
and substituting in $V_2V_6\approx V_3 V_5$ we get $V_4/V_2\approx (V_3/V_2)^2$. Continuing in this way 
yields
\begin{equation}\label{odin}
     {V_j\over V_2}\ \approx\ \left({V_3\over V_2}\right)^{j-2}\, ~{\rm for~} ~2\leq j\leq 10.
\end{equation}
From this and $F_j^2\approx V_2V_{j-2}$, we obtain  
\begin{equation}\label{dva}
    {F_j\over V_2}\ \approx\ \left({V_{3}\over V_2}\right)^{j/2-2}\, ~{\rm for~} ~5\leq j\leq 12.
\end{equation}
Substituting (\ref{odin}) and (\ref{dva}) in $f_j=V_j+c_jF_{2j}$ (recall (\ref{V})) yields
\begin{equation}\label{tri}
    {f_j\over V_2}\ \approx\ \left({V_{3}\over V_2}\right)^{j/2-2}\, ~{\rm for~} ~3\leq j\leq 6.
\end{equation}
\par This offers two ways of expressing $F_j$ for $j=5,6$: using (\ref{dva}) or (\ref{tri}), 
but with the argument $2t$ for the latter. The first gives
$$
   F_5\ \approx\ V_2(t)\left({V_3(t)\over V_2(t)}\right)^{1/2}\,,
            ~~~F_6\ \approx\ V_2(t)\left({V_3(t)\over V_2(t)}\right),
$$
and the second gives
$$
   F_5\ \approx\ V_2(2t)\left({V_3(2t)\over V_2(2t)}\right)^{3}\,,
           ~~~F_6\ \approx\ V_2(2t)\left({V_3(2t)\over V_2(2t)}\right)^4.
$$
It follows that 
$$
    {F_6\over F_5}\ \approx\ \left({V_3(t)\over V_2(t)}\right)^{1/2}\ \approx\ {V_3(2t)\over V_2(2t)}.
$$
Applying Lemma \ref{sqrtlemma} to $g(t)=V_3(t)/V_2(t)$ shows that this must converge, hence
from (\ref{tri}) the ratio $\Phi_4(t)/\Phi_3(t)$ must converge too.
Parts (i), (ii), (iii) of the theorem now follow from Lemma \ref{ratio1}, and part (iv) follows 
by de-Poissonization.
\ep

%\medskip
Combining Theorem \ref{forceRV} and Lemma \ref{noRV} we arrive at a very simple test  for the convergence,
which is easy to check in the examples of Section 4:

\begin{corollary} The condition $\lim_{j\to\infty} p_{j+1}/p_j=1$ is necessary for 
the convergence of the (normalized and centred) $X_n$ to a multivariate normal law.
\end{corollary}
\bigskip
\noindent
It should be stressed that 
the condition is by no means sufficient. For instance, the frequencies $p_j=c\{2+\sin(\log j)\}/j^2$ satisfy $p_{j+1}/p_j\to 1$ but do not
have the property of regular variation due to the oscillating sine factor. Thus in this case $X_n$ has no distributional limit.

\bigskip
\noindent
{\bf Acknowledgement}
The authors would like to thank Adrian R\"ollin and Bero Roos for helpful
discussions.  ADB gratefully acknowledges financial support from Schweizerischer 
Nationalfonds Projekt Nr.~20-117625/1.

\end {document}